\documentclass[11pt,citesort,rotate,psfig]{article}  %  ?{article}
\usepackage{epsfig,amsfonts}
\catcode`\"=13      % aktivira znak "
\def"#1{\v{#1}}      % sedaj je "c lahko \v{c}

\def\be{\begin{equation}}
\def\ee{\end{equation}}
\def\bea{\begin{eqnarray}}
\def\eea{\end{eqnarray}}

\def\ppp#1 {#1^{\prime \prime \prime}}

\def\hfl#1#2{\smash{\mathop{\hbox to 12 mm{\rightarrowfill}} \limits^{\scriptstyle #1}_{\scriptstyle #2}}}

\def\binom#1#2{\left ( {{#1} \atop {#2}} \right )}
\def \eqalign#1{\null\,\vcenter{\openup\jot \ialign{\strut\hfill$\displaystyle{##}$& $\displaystyle{{}##}$\hfill\crcr #1\crcr}}\,}
\newcommand {\qed} {\null \hfill \rule {2.5mm}{2.5mm}}
\newcommand {\R} {\ensuremath{\mathbb{R}}}
\newcommand {\N} {\ensuremath{\mathbb{N}}}

\author{Tomislav Do\v{s}li\'c$^{\star\dagger}$ and Darko Veljan$^{\ddagger}$ }
\title{Calculus proofs of some combinatorial inequalities}
\date{\today}

\begin{document}
\maketitle
$\dagger$       Department of Informatics and Mathematics, Faculty of Agriculture, University of Zagreb, Sveto\v simunska c. 25, Zagreb, CROATIA
\\ \\
$\ddagger$        Department of Mathematics , University of Zagreb, Bijeni\v cka 30, Zagreb, CROATIA\\
\\
      $\star$ To whom correspondence should be addressed,
        e-mail : doslic@faust.irb.hr\\
\newpage
{\bf \Large Abstract}
\\ \\

Using calculus we show how to prove some  combinatorial inequalities of the
type log-concavity or log-convexity. It is shown by this method that binomial
coefficients and Stirling numbers of the first and second kinds are 
log-concave, and that Motzkin numbers and secondary structure numbers of rank
$1$ are log-convex. In fact, we prove via calculus a much stronger result
that a natural continuous ``patchwork'' (i.e. corresponding dynamical systems) of
Motzkin numbers and secondary structures recursions are increasing functions. 
We indicate how to prove asymptotically the log-convexity for general secondary
structures. Our method also applies to show that 
sequences of values of some orthogonal polynomials, and in particular the
sequence of central Delannoy numbers, are log-convex.



\newpage
%------------------------------------------------------------
\section{Introduction}
%------------------------------------------------------------

In combinatorics the most prominent question is usually to find explicitly the size of certain finite set defined in an intricate way. It often happens that there
is no explicit expression for the size in question, but instead one can find
recursion, generating function or other gadgets which enable us to compute 
concrete sizes or numbers. The next question then usually asks how the sequence
of numbers satisfying certain recursion behaves. By behavior of the sequence
$\left ( a_n \right )_{n \geq 0}$ of positive real numbers it is often meant
its log-concavity (or log-convexity). Recall that a sequence 
$\left ( a_n \right )_{n \geq 0}$ of positive real numbers is {\bf log-concave}
if $a_n^2 \geq a_{n-1} a_{n+1}$ for all $n\geq 1$, and {\bf log-convex} if
$a_n^2 \leq a_{n-1} a_{n+1}$ for all $n\geq 1$. We say that a sequence
$\left ( a_n \right )_{n \geq 0}$ is {\bf log-straight} or {\bf geometric} if
$a_n^2 = a_{n-1} a_{n+1}$ for all $n\geq 1$. A (finite) sequence of positive
numbers $a_0, a_1, \ldots ,a_n$ is said to be {\bf unimodal} if, for some
$0 \leq j \leq n$ we have $a_0 \leq a_1 \leq \ldots \leq a_j \geq a_{j+1} \geq
\ldots \geq a_n$. This place $j$ is called a {\bf peak} of the sequence if it
is unique. If there are more such maximal values, we speak about a {\bf plateau} of
the sequence. It is easy to see that a log-concave positive sequence is
unimodal. The literature on log-concavity and unimodality is vast. We refer the
interested reader to the book \cite{karlin}. Combinatorial inequalities, and
in particular, the questions concerning log-concavity (or log-convexity) are
surveyed in \cite{brenti}, \cite{stanley00} and \cite{stanley89}. Some analytic 
methods are described in \cite{bender}.

In combinatorics, a preferable way to prove a combinatorial inequality is to
give a {\bf combinatorial proof}. There are two basic ways to do it. Suppose
that we are given finite sets $A$ and $B$ with $|A|=a$ and $|B|=b$ and we want
to prove, say, $a \leq b$. One way to prove it is to construct an injection
$A \rightarrow B$ (or a surjection $B \rightarrow A$), and the other is to
show that the number $c=b-a$ is nonnegative, by showing that $c$ is cardinality
of certain set or that $c$ is the dimension of certain vector space (and hence
nonnegative) etc. As an example, let us show that binomial coefficients $
\binom {n}{k}$, $k=0,1,\ldots ,n$ are log-concave. It is trivial to check
algebraically that $\binom {n}{k} ^2 \geq \binom {n}{k-1} \binom {n}{k+1}$ by
using the standard formula $\binom {n}{k} = \frac {n!}{k!(n-k)!}$, but 
combinatorially it goes as follows.

First define the {\bf Narayana numbers} $N(n,k)$ for integers $n, \quad k \geq
1$ as
$$N(n,k)=\frac {1}{n} \binom {n}{k} \binom {n}{k-1} = \frac {1}{k}\binom {n}{k-1}\binom{n-1}{k-1},$$
and $N(0,0):=1$. Next we note that 
$$\binom {n}{k} ^2-\binom {n}{k-1} \binom {n}{k+1} = \left |
{{\binom {n}{k} \atop \binom {n}{k-1}} \quad {\binom {n}{k+1} \atop \binom {n}{k}}} \right | = N(n+1,k+1).$$
Finally, we need the fact that Narayana numbers have a combinatorial meaning, i.e.
they count certain finite sets (see below). Therefore we get $\binom {n}{k} ^2
- \binom {n}{k-1} \binom {n}{k+1} \geq 0.$ There are also other combinatorial 
proofs of log-concavity of binomial coefficients, as well as log-concavity of
Stirling numbers (of both kinds) etc., but they are all rather involved and/or
tricky. In this paper we present a way to prove various combinatorial inequalities
by a straightforward method of calculus. Inductive and injective proofs of
log-convexity results are described in \cite{sagan}.

\section {Calculus proofs of log-concavity and log-convexity properties}

Let us first recall briefly calculus proofs of log-concavity of binomial
coefficients and Stirling numbers. Let $c(n,k)$ be the number of permutations 
of the set $[n]:=\{1,2,\ldots ,n\}$ with exactly $k$ cycles and $S(n,k)$ the
number of partitions of $[n]$ into exactly $k$ parts (or blocks). The numbers
$c(n,k)$ and $S(n,k)$ are called {\bf Stirling numbers} of the {\bf first} and
{\bf second kind}, respectively. The following formulae are well known (see 
\cite{stanley}).
\be
(x+1)^n = \sum_{k=0}^n \binom {n}{k} x^k,
\ee
\be
x^{\bar n}=x(x+1)\ldots (x+n-1) = \sum_{k=0}^n c(n,k) x^k,
\ee
\be
x^n= \sum_{k=0}^n S(n,k) x^{\underline k},
\ee
where $x^{\underline k}:=x(x-1)\ldots (x-k+1)$ is the $k$-th {\bf falling power}
and $x^{\bar k} = x(x+1)\ldots (x+k-1)$ the $k$-th {\bf rising power} of x.

The following Newton's lemma is a consequence of the Rolle's theorem from
calculus.

{\bf Lemma 1.}\\ Let $P(x)=\sum_{k=0}^n a_k x^k$ be a real polynomial whose
all roots are real numbers. Then its coefficients are log-concave, i.e.
$a_k^2 \geq a_{k-1} a_{k+1}, k=1, \ldots ,n-1$. (Moreover, $\frac {a_k}{\binom
{n}{k}} $ are log-concave). \qed

Now, from (1) and (2) we see that $(x+1)^n$ and $x^{\bar n}$ have only real roots
and by Lemma 1. we conclude that the sequences $\binom {n}{k}$ and $c(n,k)$ are
log-concave.

The case of the sequence $S(n,k)$ is a bit more involved. We claim that the
polynomial
\be 
P_n(x)=\sum_{k=0}^n S(n,k) x^k
\ee
has all real roots (in fact non-positive and different). Namely, $P_0(x)=1$
and from the basic recursion $$S(n,k)=S(n-1,k-1)+kS(n-1,k)$$ it follows at once
that $$P_n(x)=x\left [P_{n-1}'(x)+P_{n-1}(x)\right ].$$
The function $Q_n(x)=P_n(x)e^x$ has the same roots as $P_n(x)$ and it is easy to
verify $Q_n(x)=xQ_n'(x)$. By induction on $n$ and by using the Rolle's theorem
it follows easily that $Q_n$ and hence $P_n$ have only real and non-positive
roots.

So, we have proved by calculus the following.

{\bf Theorem 1.}\\The sequences $\binom {n}{k}_{k \geq 0}, \left ( c(n,k) \right)_{k \geq 0}, \left (S(n,k) \right )_{k \geq 0}$ are log-concave. Hence they
are also unimodal. \qed

It is also well known that the peak of the sequence $\binom {n}{k}$ is at
$k=\lfloor n/2 \rfloor$, while the peak for the other two sequences is much 
harder to determine. It is known that $S(n,k)$'s reach their peak for $k \approx
n/\log n$, if $n$ is large enough. (An inductive proof of Theorem 1. is given in
\cite{sagan}.)

Now we turn to a different kind of combinatorial entities. Recall that a {\bf 
Dyck path} is a path in the coordinate $(x,y)$-plane from$(0,0)$ to $(2n,0)$
with steps $(1,1)$ and $(1,-1)$ never falling below the $x$-axis. Denote the set
of all such paths by ${\cal D}_n$. A {\bf peak} of a path $P \in {\cal D}_n$ is
a place at which the step $(1,1)$ is directly followed by the step $(1,-1)$.
Denote by ${\cal D}_{n,k} \subseteq {\cal D}_n$ the set of all Dyck paths of
length $2n$ with exactly $k$ peaks. Note that $1 \leq k \leq n$. The following
facts are also well known (see \cite{stanley}).
$$|{\cal D}_n| = \frac {1}{n+1} \binom {2n}{n} = C_n $$
$$|{\cal D}_{n,k}| = N(n,k),$$
where $C_n$ is $n$-th {\bf Catalan number}. The Catalan numbers are log-convex.
The Narayana numbers are log-concave in $k$ for fixed $n$. Both these facts can easily
be proved algebraically, but there are also combinatorial proofs, as well as
calculus proofs. We omit here these proofs, since we want to emphasize the 
following more intricate combinatorial quantities, related to the above just
introduced. 

A {\bf Motzkin path} is a path in the coordinate $(x,y)$-plane from $(0,0)$ to $(n,0)$
with steps $(1,1)$, $(1,0)$ and $(1,-1)$ never falling below the $x$-axis. Let
${\cal M}_n$ be the set of all such paths and let $M_n = |{\cal M}_n|$. The 
number $M_n$ is called the $n$-th {\bf Motzkin number}. 

Some basic properties of Motzkin numbers are as follows (\cite{stanley},
 \cite{doslic1}).

{\bf Theorem 2.}\\
(a) $M_n=\sum_{k \geq 0} \binom {n}{2k}C_k, \quad C_{n+1}=\sum_{k \geq 0}\binom {n}{k} M_k;$\\
(b) $M_{n+1}=M_n+\sum_{k=0}^{n-1} M_kM_{n-k-1};$\\
(c) The generating function of $\left ( M_n \right )_{n\geq 0}$ is given by
$$M(x)=\sum_{n\geq 0}M_n x^n = \frac {1-x-\sqrt {1-2x-3x^2}}{2x^2};$$
(d) $(n+2)M_n=(2n+1)M_{n-1}+3(n-1)M_{n-2};$\\
(e) $M_n \sim \sqrt {\frac {3}{4 \pi }} 3^{n+1} n^{-3/2}.$ \qed

The log-convexity of the sequence of Motzkin numbers was first established
algebraically in \cite{aigner}, and shortly afterwards combinatorial proof
appeared in \cite{callan}.
We shall prove now by calculus that $\left ( M_n \right )_{n\geq 0}$ is a 
log-convex sequence and some consequences of this property.

{\bf Theorem 3.}\\
(a) The sequence $\left ( M_n \right )_{n\geq 0}$ is log-convex;\\
(b) $M_n \leq 3 M_{n-1}$, for all $n\geq 1$;\\
(c) There exists $x=\lim_{n \rightarrow \infty} \frac {M_n}{M_{n-1}}$, and $x=3$.

{\bf Proof}\\
(a) Let us start from the short recursion in Theorem 2.(d):
$$M_n=\frac {2n+1}{n+2}M_{n-1} +\frac {3(n-1)}{n+2}M_{n-2}.$$
Divide this recursion by $M_{n-1}$ and denote $x_n:= \frac {M_n}{M_{n-1}}$. Then
we obtain the following recursion:
\be 
x_n=\frac {2n+1}{n+2} +\frac{3(n-1)}{n+2}\frac {1}{x_{n-1}}
\ee
with initial condition $x_1=1$. The log-convexity $M_n^2 \leq M_{n-1}M_{n+1}$
is equivalent to $x_n \leq x_{n+1}$. To prove that $(x_n)_{n\geq 0}$ is an increasing
sequence, we shall prove a much stronger claim. To this end, define the following
function $f : [2,\infty) \rightarrow \R$. For $x\in [2,3]$, define $f(x)=2$.
For $x\geq 3$, let (by simulating (5))
\be
f(x)=\frac {2x+1}{x+2} +\frac {3(x-1)}{x+2} \frac {1}{f(x-1)}.
\ee
Note that $f(n)=x_n$. We shall prove that $f$ is an increasing function, and
consequently that $\left ( x_n \right )_{n\geq 0}$ is an increasing sequence.
Note first that the function $f$ is continuous ($f$ is, in fact, a dynamical
system), and on every open interval $(n,n+1)$, where $n\geq 2$ is an integer,
$f$ is a rational function, with no poles on it. Therefore, $f$ is smooth on 
every open interval $(n,n+1)$, for $n\geq 2$. Note that, for example, $f(x)=
\frac {7x-1}{2(x+2)}$ for $x\in [3,4]$, $f(x)=\frac {20x^2-9x-14}{7x^2+6x-16}$
for $x\in [4,5]$, etc. It is trivial to check that $f(x) \geq 2$, for all 
$x\geq 2$. Suppose inductively that $f$ is an increasing function on a segment
$[3,n]$. For $n=4$ it is (almost evidently) true. Let $n \geq 4$, and take a
point $x\in (n,n+1)$. By taking the derivative $f'(x)$ of (6), and plugging in
once more the term for $f'(x-1)$, we have:
$$\displaylines {
\qquad f'(x)=\frac {3}{[(x+2)f(x-1)]^2} \Bigg [f^2(x-1)+3f(x-1)- 3 \frac{(x-1)(x+2)} {(x+1)^2f(x-2)}[f(x-2)+3] \hfill \cr 
\hfill +3 \frac{(x^2-1)(x^2-4)}{[(x+1)f(x-2)]^2} f'(x-2) \Bigg ] \qquad \cr }$$

By inductive hypothesis, $f$ is an increasing function on $[3,n]$ and hence
$f(x-1) \geq f(x-2) \geq 2$ and $f'(x-2) \geq 0$. So, it is enough to prove
that $f'(x)\geq 0$. However, this follows from the following.

The last term in square brackets is clearly positive, by the induction 
hypothesis. We claim that the rest is positive, too. This claim is equivalent
with
$$\frac {[f^2(x-1)+3f(x-1)]f(x-2)}{f(x-2)+3} \geq 3 \frac {(x-1)(x+2)}{(x+1)^2}.$$
But this inequality is true, since by inductive hypothesis $f(x-1) \geq f(x-2)
\geq 2$, and hence the left hand side is at least equal to $f^2(x-2) \geq 4$,
while the right hand side has the maximum (for $x \geq 3$) equal to $3$.
Hence $f'(x) > 0$ for all $x \in (n,n+1)$. So, the function $f$ is
strictly increasing on $(n,n+1)$, and then, by continuity, also on $[3,n+1]$.
In particular, $f(n+1)=x_{n+1} \geq x_n = f(n)$.
This completes the step of induction.

(b) and (c) follow now simultaneously, because by (a), the sequence 
$\left ( x_n \right )_{n\geq 0}$ is increasing and from (5) it follows easily
by induction on $n$ that $2 \leq x_n \leq 7/2$, i.e. $\left ( x_n \right )$ is
bounded. \qed

Closely related combinatorial structures to Motzkin paths are the so called
secondary structures. A {\bf secondary structure} is a simple planar graph
on vertex set $[n]$ with two kinds of edges: segments $[i,i+1]$, for $1\leq i
\leq n-1$ and arcs in the upper half-plane which connect some $i,j$, where
$i < j$ and $j-i > l$, for some fixed integer $l\geq -1$, such that the arcs 
are totally disjoint. Such a structure
is called a secondary structure of {\bf size} $n$ and {\bf rank} $l$. The
importance for the study of these structures comes from biology. They are 
crucial in understanding the role of RNA in the cell metabolism and in decoding
the hereditary information contained in DNA. Biologists call the vertices of a 
secondary structure {\bf bases}, the segments they call {\bf p-bonds} (p stands
for phosphorus) and arcs they call {\bf h-bonds} (h stands for hydrogen). Let
${\cal S}^{(l)}(n)$ be the set of all secondary structures of rank $l$ on $n$
vertices and $S^{(l)}(n)=|{\cal S}^{(l)}(n)|$ the {\bf secondary structure numbers} of
rank $l$. In a sense, the Motzkin numbers are secondary structure numbers of rank $0$,
and the Catalan numbers are secondary structure numbers of the (degenerate) rank $-1$.
In these cases the corresponding graphs are not simple, but the other requirements
on secondary structures remain.

Now we shall apply our method of calculus to prove that in the case 
$l=1$ the behavior of the numbers $S^{(1)}(n)$ is also log-convex. So, we have:

{\bf Theorem 4.}\\ The sequence $\left (S^{(1)}(n) \right )_{n\geq 0}$ is
log-convex.

{\bf Proof}\\ As for the Motzkin numbers, it turns out that for $S^{(1)}(n)$ the
following short recursion holds (see \cite{doslic1} and \cite{doslic}):
\be
(n+2)S^{(1)}(n)=(2n+1)S^{(1)}(n-1)+(n-1)S^{(1)}(n-2)+(2n-5)S^{(1)}(n-3)-(n-4)
S^{(1)}(n-4)
\ee
with initial conditions $S^{(1)}(0)=S^{(1)}(1)=S^{(1)}(2)=1, S^{(1)}(3)=2$.
By dividing this recursion with $S^{(1)}(n-1)$ and denoting
$$x_n=\frac {S^{(1)}(n)}{S^{(1)}(n-1)},$$
we get
\be
x_n=\frac {1}{n+2} \left [2n+1 + \frac {n-1}{x_{n-1}} +\frac {2n-5}{x_{n-1}x_{n-2}}-\frac {n-4}{x_{n-1}x_{n-2}x_{n-3}}\right ],
\ee
with initial conditions $x_3=x_4=x_5=2$ (note that $x_1=x_2=1$).

The log-convexity of $S^{(1)}(n)$'s is equivalent with the fact that $(x_n)$
is an increasing sequence.

Now define the function $f:[2,\infty ) \rightarrow \R$ by simulating (8) as:
\be
f(x)= \cases {
2 &, if $x\in [2,5]$,\cr
\frac {1}{x+2}\left [2x+1 + \frac {x-1}{f(x-1)} +\frac {2x-5}{f(x-1)f(x-2)}-
\frac {x-4}{f(x-1)f(x-2)f(x-3)}\right ] &, if $x\geq 5.$\cr
}
\ee

Clearly, for any integer $n\geq 3$, $f(n)=x_n$, and $f$ is continuous, and,
in fact, piecewise rational and smooth on any open interval $(n,n+1)$ for
$n\geq 2$. The basic idea is, as in the proof of Theorem 3.(a), to show that
$f$ is an increasing and bounded function, and hence $(x_n)$ is an increasing sequence.
In next few lemmas we proceed with details.

{\bf Lemma 2.}\\ For all $x\geq 2$, we have $2 \leq f(x) \leq 3$, while for
$x\geq 53$ we have even stronger bounds:\\ $2.5 \leq f(x) \leq 2.67$.

{\bf Proof} \\ We prove inductively that $2 \leq f(x) \leq 3$ for $x\in [2,n]$.
For $n\leq 11$ it can be checked directly. Let $n\geq 11$ and $x\in (n,n+1]$. Then
$$\displaylines{
\qquad f(x) \leq \frac {1}{x+2}\left [2x+1 + \frac {x-1}{f(x-1)} +\frac {2x-5}{f(x-1)f(x-2)}\right ] \hfill \cr
\hfill \leq \frac {1}{x+2}\left [2x+1 + \frac {x-1}{2}+\frac {2x-5}{4}\right ]
=\frac {12x-3}{4x+8} \leq 3. \qquad\cr
}$$
On the other hand, $$f(x) \geq \frac {1}{x+2}\left [2x+1 + \frac {x-1}{3}+\frac {2x-5}{9}-\frac {x-4}{8} \right ]
=\frac {175x+44}{72x+144} \geq 2,$$
for all $x\geq 8$. So, $2 \leq f(x) \leq 3$ on $(n,n+1]$ and the first claim is
proved.

The stronger bounds also follow by induction. By direct computation, (using
Mathematica) one can check that they hold on the interval $[53,56]$. Suppose
$2.5 \leq f(x) \leq 2.67$ on some interval $[53,n]$, where $n\geq 56$ and take
$x\in[n,n+1]$. From (8) we get
$$f(x) \leq \frac {1}{x+2}\left [2x+1 + \frac {x-1}{2.5}+\frac {2x-5}{2.5^2}
-\frac {x-4}{2.67^3}\right ] = \frac {2.6675x+0.010148}{x+2} \leq 2.67,$$
for all $x\geq 0$. On the other hand,
$$f(x) \geq \frac {1}{x+2}\left [2x+1 + \frac {x-1}{2.67}+\frac {2x-5}{2.67^2}
-\frac {x-4}{2.5^3}\right ] = \frac {2.59108x+0.0181}{x+2},$$
and this is greater than $2.5$ for $x \geq 53$ (since the right hand side is
equal to $2.5$ for $x=52.918$). So, Lemma 2. is proved. \qed

{\bf Lemma 3.}\\ The function $f$ is increasing.

{\bf Proof}\\ Suppose again inductively on $n \in \N$ that $f$ increases on $
(5,n]$. We shall prove that $f$ increases on $(n,n+1)$. One can check directly
(using, e.g. Mathematica) that $f$ increases on $(5,n_0]$, as far as $n_0=61$.
Namely, the function $f$ on interval $(n,n+1)$ is a rational function whose
both numerator and denominator are polynomials with integer coefficients of
degree $n-4$. The derivative of $f$ is also a rational function, and its 
denominator is always positive. So, we need to show that the numerators of the
derivative of $f$ are positive on every interval $(n,n+1)$, for $n\leq n_0-1$.
An advanced computer algebra system, such as {\it Mathematica}, gives us
readily explicit expressions for $f(x)$ and $f'(x)$ on any given interval
$(n,n+1)$. Let us denote $f'(x)=\frac {N_n(x)}{D_n(x)}$ on interval $(n,n+1)$.
If we can find some $k \in \N$, $k\leq n$, such that all coefficients of 
$N_n(x+k)$ are nonnegative, we are done, since then $f'(x)$ can not change its
sign on the considered interval. It turns out that $k=2$ works for all intervals
$(n,n+1)$ with $n \leq 60$. Hence, $f'(x)\geq 0$ for $x\in (n,n+1)$, $n \leq 60$
and $f(x)$ is increasing on $[5,61]$. It is important to note here that all 
performed computations include only integer quantities, and no round-off errors
occur. 

Take $x\in (n,n+1)$ for $n\geq n_0$. Then $f'(x)>0$ for $x\in (i,i+1), \quad
i=5, \ldots ,n-1$, and also $f(x)\geq f(x-1)$, for $4 \leq x \leq n$.

Denote for short $f_i=f(x-i)$, $i\geq 1$. Then (9) can be written as
$$(x+2)f_1f_2f_3f(x)=(2x+1)f_1f_2f_3+(x-1)f_2f_3+(2x-5)f_3-(x-4).$$
By taking derivative, we get
$$f'(x)=\frac {1}{D(x)}\left [F(x)+F_3(x)f_3'(x)-F_1(x)f_1'(x)-F_2(x)f_2'(x)\right ],$$
where 
$$\eqalign{
D(x)& =(x+2)f_1f_2f_3,\cr
F(x) &=2f_1f_2f_3+f_2f_3+2f_3-1-f_1f_2f_3f(x),\cr
F_1(x) &=\left [(x+2)f(x)-(2x+1)\right ]f_2f_3,\cr
F_2(x) &=\left [(x+2)f_1f(x)-(2x+1)f_1 -(x-1)\right ]f_3,\cr
F_3(x) &=2x-5+(2x+1)f_1f_2+(x-1)f_2-(x+2)f_1f_2f(x).\cr
}$$

Using (9), let us express $D(x), F(x), F_i(x), i=1,2,3$ only in terms of $f_i$'s
and $x$: \\
$$\eqalign{
D(x) & =(x+2)f_1f_2f_3,\cr
F(x) & =\frac {3}{x+2}(f_1f_2f_3+f_2f_3+3f_2-2),\cr
F_1(x) & =\frac {1}{f_1}\left [(f_2f_3+2f_3-1)x-(f_2f_3+5f_3-4)\right ],\cr
F_2(x) & =\frac {1}{f_2}\left [(2f_3-1)x-(5f_3-4)\right ],\cr
F_3(x) & =\frac {1}{f_3} (x-4).\cr
}$$

Now plug in derivatives $f_1'$ and $f_2'$ by the same rule, to obtain
%$$\displaylines{
%\qquad f'(x)=\frac {1}{D(x)}\Bigg \{F(x) - \frac {F_1(x)F(x-1)}{D(x-1)}-
%\frac {F_2(x)F(x-2)}{D(x-2)}+\frac {F_1(x)F_1(x-1)}{D(x-1)}f_2'+ \qquad\cr
%\hfill \left [ \frac {F_1(x)F_2(x-1)}{D(x-1)}+\frac {F_2(x)F_1(x-2)}{D(x-2)}+F_3(x)\right ]f_3'+\hfill (10)\cr
%\hfill \left [\frac {F_2(x)F_2(x-2)}{D(x-2)}-\frac {F_1(x)F_3(x-1)}{D(x-1)}\right ]f_4'
%-\frac {F_2(x)F_3(x-2)}{D(x-2)}f_5' \Bigg \}\qquad\cr
%}$$
\bea
f'(x) & = & \frac {1}{D(x)}\Bigg \{F(x) - \frac {F_1(x)F(x-1)}{D(x-1)}-
\frac {F_2(x)F(x-2)}{D(x-2)}+\frac {F_1(x)F_1(x-1)}{D(x-1)}f_2' \nonumber \\
 & + & \left [ \frac {F_1(x)F_2(x-1)}{D(x-1)}+\frac {F_2(x)F_1(x-2)}{D(x-2)}+F_3(x)\right ]f_3' \cr
 & + & \left [\frac {F_2(x)F_2(x-2)}{D(x-2)}-\frac {F_1(x)F_3(x-1)}{D(x-1)}\right ]f_4'
-\frac {F_2(x)F_3(x-2)}{D(x-2)}f_5' \Bigg \}.
\eea
The ``coefficients'' by $f_2'$ and $f_3'$ are positive. By further pumping in
$f_5'$, the terms $f_6'$ and $f_7'$ will appear with positive ``coefficients'',
while $f_8'$ will appear with negative ``coefficient'' and a ``free'' negative
term $$-\frac {F_2(x)F_3(x-2)F(x-5)}{D(x-2)D(x-5)}$$ also appears. Every further 
pumping in for $f_{3k+2}'$ contributes positive terms by $f_{3k+3}'$ and 
$f_{3k+4}'$, a negative term by $f_{3k+5}'$ and a negative ``free'' term. If
we continue to pump in long enough, the argument of the negative term will
be eventually ``trapped'' in the interval $(2,5)$, and there $f'=0$. So, to
prove that $f'(x)>0$ we only have to show that the ``coefficient'' of $f_4'$
is positive and that ``free'' term (i.e. the term without any $f_i'$) is also
positive. These two facts we prove in the next lemma.

{\bf Lemma 4.}\\The ``coefficient'' of $f_4'$ and the ``free'' term, obtained
by pumping in $f_5', f_8', \ldots$ in (10) are both positive. More precisely,
with previous notations we have:\\
$(a)$ \hskip 5cm $
L_4(x):=\frac {F_2(x)F_2(x-2)}{D(x-2)}-\frac {F_1(x)F_3(x-1)}{D(x-1)} \geq 0,
$\\
for $x\geq n_0$;
$$\displaylines{
(b) \quad\quad
\qquad L(x):=F(x)-\frac {F_1(x)F(x-1)}{D(x-1)}-\frac {F_2(x)F(x-2)}{D(x-2)} \hfill \cr
\hfill -\frac {F_2(x)F_3(x-2)}{D(x-2)} \frac{F(x-5)}{D(x-5)}\left [1+
\frac{F(x-8)}{D(x-8)}+\frac{F(x-8)}{D(x-8)}\frac{F(x-11)}{D(x-11)}+\ldots \right] \geq 0,\qquad\cr
}$$
for $x\geq n_0$, where $n_0$ can be taken in the worst case to be $n_0=61$.

{\bf Proof}\\ (a) The condition $L_4(x) \geq 0$ is easily seen to be equivalent
to
$$\displaylines{
\qquad (x+1)f_1\left [(2f_3-1)x-(5f_3-4)\right ]\left [(2f_5-1)(x-2)-(5f_5-4)\right ]\hfill\cr
\hfill -xf_5\left [(f_2f_3+2f_3-1)x-(f_2f_3+5f_3-4)\right ](x-5) \geq 0.\qquad\cr
}$$
If we leave out the factor $(x+1)f_1$ from the first term and the factor $xf_5$
from the second term, we obtain even stronger inequality (recall, we are still
under inductive hypothesis, and this implies that $f_1 \geq f_5$). By grouping
terms by powers of $x$, this stronger inequality can be written in the form
$$c_{24}(x)x^2+c_{14}(x)x+c_{04}(x)=[c_{24}(x)x+c_{14}(x)]x+c_{04}(x)\geq 0,$$
where \\
$$\eqalign{
c_{24}(x) & =4f_3f_5-f_2f_3-4f_3-2f_5+2,\cr
c_{14}(x) & =6f_2f_3+17f_5+32f_3-28f_3f_5-19,\cr
c_{04}(x) & =45f_3f_5-5f_2f_3-36f_5-55f_3+44.\cr
}$$
Now estimate $c_{24}(x), \quad c_{14}(x)$ and $c_{04}(x)$ using the bounds from 
Lemma 2. We easily obtain $c_{24}(x) \geq 3.8516, \quad c_{14}(x) \geq -58.6042
\quad c_{04}(x) \geq 46.6355$ for $x\geq n_0$. For example, since $f_{min}=2.5$,
$f_{max}=2.67$ for $x\geq n_0$, we have then
$$c_{24}(x) \geq 4f_{min}^2-f_{max}^2-6f_{max}+2 = 3.8516.$$
These bounds then imply $c_{24}(x)x+c_{14}(x) \geq 0$, and hence 
$[c_{24}(x)x+c_{14}(x)]x+c_{04}(x)\geq 0$, for $x\geq n_0$. So, $L_4(x) \geq 0$
for $x\geq n_0$ and the claim (a) is proved.

(b) First of all, the function $\frac {F(x)}{D(x)}$ is easily seen to be less
than $\frac {129}{8(x+2)^2}$ (by using $2 \leq f_1, f_2, f_3 \leq 3$). For
$x \geq 10$, it follows then that $$\frac {F(x)}{D(x)} \leq q,$$ where $q =
\frac {129}{1152}$. By using $\frac {F(x-i)}{D(x-i)}\leq q$ in the brackets
of (b), we see that this sum is less than the sum of the geometric series
$1+q+q^2+\ldots =\frac {1}{1-q} < 2$. Hence $L(x) \geq 0$ will be a consequence
of the stronger inequality:
$$F(x)-\frac {F_1(x)F(x-1)}{D(x-1)}-\frac {F_2(x)F(x-2)}{D(x-2)}-2\frac {F_2(x)F_3(x-2)}{D(x-2)} \frac{F(x-5)}{D(x-5)}\geq 0.$$
But, since we do not know which one of the quotients $\frac {F(x-1)}{D(x-1)}$,
$\frac {F(x-2)}{D(x-2)}$ and $\frac {F(x-5)}{D(x-5)}$ is the largest, the last
inequality will be a consequence of the three inequalities in the next Lemma.

\newpage
{\bf Lemma 5.}\\ Keeping the same notations as above, we have
$$\displaylines{
(a) \quad\quad
F(x)\geq \left [F_1(x)+F_2(x)+2\frac {F_2(x)F_3(x-2)}{D(x-2)}\right ]\frac{F(x-5)}{D(x-5)}, \quad x\geq n_0,\qquad\qquad\qquad\qquad\qquad\qquad\qquad\cr
(b) \quad\quad
F(x)\geq \left [F_1(x)+F_2(x)+2\frac {F_2(x)F_3(x-2)}{D(x-2)}\right ]\frac{F(x-2)}{D(x-2)}, \quad x\geq n_0,\qquad\qquad\qquad\qquad\qquad\qquad\qquad\cr
(c) \quad\quad
F(x)\geq \left [F_1(x)+F_2(x)+2\frac {F_2(x)F_3(x-2)}{D(x-2)}\right ]\frac{F(x-1)}{D(x-1)}, \quad x\geq n_0.\qquad\qquad\qquad\qquad\qquad\qquad\qquad\cr
}$$

{\bf Proof}\\ We shall prove only (a) with substantial details. The other two
inequalities can be proved essentially in the same manner. The inequality (a)
is equivalent to
$$x(x-3)^2f_1f_2f_3f_4f_5^2f_6f_7f_8A\geq [xf_3f_4f_5^2(ax-b)+2f_1(cx-d)(x-6)](x+2)B,$$
where \\
$$\eqalign{
A & =f_1f_2f_3+f_2f_3+3f_3-2,\cr
B & =f_6f_7f_8+f_7f_8+3f_8-2,\cr
a & =f_2^2f_3+2f_1f_3+2f_2f_3-f_1-f_2,\cr
b & =f_2^2f_3+5f_1f_3+5f_2f_3-4f_1-4f_2,\cr
c & =2f_3-1, \cr
d & =5f_3-4.
}$$
By inductive hypothesis it follows that $A\geq B$, and so if we prove the stronger 
inequality by leaving out $A$ and $B$ in the above inequality, we are done.
But this stronger inequality turns out to be (after grouping terms by powers
of $x$ and some manipulations):
$$c_{35}(x)x^3+c_{25}(x)x^2+c_{15}(x)x+c_{05}(x) \geq 0,$$
or, what is the same,
\be
[c_{35}(x)x+c_{25}(x)]x^2+c_{15}(x)x+c_{05}(x) \geq 0,
\ee
where
$$c_{35}(x)=f_1f_2f_3f_4f_5^2f_6f_7f_8-f_2^2f_3^2f_4f_5^2-2f_2f_3^2f_4f_5^2-
2f_1f_3^2f_4f_5^2+f_2f_3f_4f_5^2+f_1f_3f_4f_5^2-4f_1f_3+2f_1,$$
$$c_{25}(x)=-6f_1f_2f_3f_4f_5^2f_6f_7f_8-f_2^2f_3^2f_4f_5^2+f_2f_3^2f_4f_5^2
+f_1f_3^2f_4f_5^2-2f_2f_3f_4f_5^2-2f_1f_3f_4f_5^2+26f_1f_3-16f_1,$$
$$c_{15}(x)=9f_1f_2f_3f_4f_5^2f_6f_7f_8+2f_2^2f_3^2f_4f_5^2
+10f_2f_3^2f_4f_5^2+10f_1f_3^2f_4f_5^2-8f_2f_3f_4f_5^2-8f_1f_3f_4f_5^2+8f_1f_3+8f_1,$$
$$c_{05}(x)=96f_1-120f_1f_3.$$
Now we estimate the above functions $c_{i5}(x)$ by the bounds from Lemma 2., 
$f_{min}=2.5$ and $f_{max}=2.67$ for $x\geq n_0$. We have
$$c_{35}(x)\geq f_{min}^9-f_{max}^7-4f_{max}^6+2f_{min}^5-4f_{max}^2+2f_{min}=
1569.9574,$$
and similarly $c_{25}(x)\geq -42278.4392$, $c_{15}(x)\geq 38334.7087$ and
$c_{05}(x)\geq -615.468$. This altogether then yields $[c_{35}(x)x+c_{25}(x)]
\geq 0$ for $x\geq n_0$, and this in turn implies (11) for $x\geq n_0$. Thus we
have proved (a).

As we said earlier, the inequalities (b) and (c) can be proved in the same way,
and we omit their proofs. \qed

To conclude, by lemmas 4. and 5. and induction hypothesis $f_i'\geq 0$ we have
shown that $f'(x)\geq 0$ for $x\in (n,n+1)$. By continuity of $f$ it follows
that $f$ is increasing on $(5,n+1)$, hence on $(5,n+1]$ and by induction $f$
is increasing on the whole interval $(2,\infty )$. This finally proves Theorem
4. \qed

This proof of Theorem 4., although rather involved (mostly computationally), is
conceptually quite simple, and can be considered as a calculus proof. Once 
again, our proofs of Theorems 1., 3. and 4. show the strong interference
between ``discrete'' and ``continuous'' mathematics.

We note finally that the proofs of Theorems 3 and 4 we have presented here
prove much stronger claims than actually stated in these theorems. Namely, they
show not only that sequences $(x_n)$ given by recursions (5) and (8) are 
increasing, but also that their natural continuous ``patch-works'' are 
increasing functions, too. Theorems 3 and 4 itself can be proved much simpler in 
such a way that we interlace the sequences $(x_n)$ given by recursions (5) and (8)
with an increasing sequence $a_n$, i.e. $a_n \leq x_n \leq a_{n+1}$. In the
case (5), $a_n=\frac {6n}{2n+3}$ for $n\geq 3$, and in the case (8) $a_n=\frac 
{2n\phi ^2}{2n+3}$, for $n\geq 6$, where $\phi = \frac {1 + \sqrt 5}{2}$ 
is the golden ratio.

This ``interlacing'' or ``sandwiching'' method can also be applied to prove
the log-convexity of sequences $S^{(l)}(n)$ for $l=2,3$ and $4$. The details
are rather involved and will appear elsewhere.

We are not aware of any combinatorial proofs of the log-convexity property of
the sequences $S^{(l)}(n)$.

It can be proved by geometric reasoning that the numbers $S^{(l)}(n)$
of rank $l$ secondary structures asymptotically behave as
$$S^{(l)}(n)\sim K_l \alpha _l^n n^{-3/2},$$
where $K_l$ and $\alpha _l$ are constants depending only on $l$, and $\alpha _l
\in [2,3]$ and $\alpha _l \searrow 2$ as $l \rightarrow \infty$. The constant
$\alpha _l$ is the largest real solution of $x^l(x-2)^2=1$. For instance,
$\alpha _0=3$, $\alpha _1=(3+\sqrt 5)/2$, $\alpha _2=1+\sqrt 2$, and $\alpha _3$,
$\alpha _4$, $\alpha _5$ and $\alpha _6$ can be also explicitly
computed (see \cite{doslic1}).

By taking the quotient $x^{(l)}_n=\frac {S^{(l)}(n)}{S^{(l)}(n-1)}$, we see
that 
$$x^{(l)}_n=\frac {S^{(l)}(n)}{S^{(l)}(n-1)} \sim \alpha _l \left (1-\frac {1}{n}\right )^{3/2} := a^{(l)}_n.$$
Clearly, the sequence $\left ( a^{(l)}_n \right )_{n \geq 1}$ increasingly 
tends to $\alpha _l$ as $n \rightarrow \infty$. This suggests that
$\left ( x^{(l)}_n \right )_{n \geq 1}$ should be interlaced with 
$\left ( a^{(l)}_n \right )_{n \geq 1}$, at least asymptotically.

These and many other properties of general secondary structures will appear
elsewhere \cite{doslic1}. More on the biological background of secondary structures
the reader can find in \cite{waterm} and \cite{kruskal}.

Our ``calculus method'' can be applied to many other combinatorial quantities
as well. For example, it can be proved in this way (see \cite{doslic1}) that
{\bf big Schr\"oder numbers} $r_n$ are log-convex. Recall that $r_n$ is the number 
of lattice paths from $(0,0)$ to $(n,n)$ with steps $(1,0)$, $(0,1)$ and $(1,1)$
that never rise above the line $y=x$.

As our final example, let us consider the sequence $P_n(t)$ of the values of
Legendre polynomials in some fixed real $t \geq 1$. We start from Bonnet's
recurrence (see \cite{szego}):
\be 
P_n(t)=\frac {2n-1}{n} tP_{n-1}(t) - \frac {n-1}{n} P_{n-2}(t), \quad n\geq 2,
\ee
with $P_0(t)=1$, $P_1(t)=t$.
Dividing this by $P_{n-1}(t)$ and denoting the quotient $\frac {P_n(t)}{
P_{n-1}(t)}$ by $x_n(t)$, we get the following recursion for $x_n(t)$:
\be 
x_n(t)=t \frac {2n-1}{n} - \frac {n-1}{n} \frac {1}{x_{n-1}(t)}
\ee
with initial condition $x_1(t)=t$. The log-convexity of the sequence $P_n(t)$
will follow if we show that the sequence $x_n(t)$ is increasing.

To this end we define the function $f_t(x) : [0,\infty) \rightarrow \R$ by
\be
f_t(x)= \cases {
t &, if $x\in [0,1]$,\cr
t\frac{2x-1}{x} - \frac {x-1}{x}\frac {1}{f_t(x-1)} &, if $x\geq 1$\cr
}
\ee
It is easy to show by induction on $n$ that $f_t$ is continuous and piecewise
rational function on any interval $[1,n]$. By the same method it easily follows
that $f_t$ is bounded, i.e. $1 \leq f_t(x) \leq 2t$ for all $x\geq 1$.
It is clear that $f_t(n)=x_n(t)$, for any integer $n\geq 1$.

{\bf Theorem 5.}\\ The sequence $P_n(t)$ of the values of Legendre polynomials
is log-convex for any fixed real $t \geq 1$.

{\bf Proof}\\ The claim will follow if we show that $f_t(x)$ is an increasing
function on $[1,\infty)$. From piecewise rationality and boundedness of $f_t$
it follows that $f_t$ is differentiable on every open interval $(n,n+1)$.
Suppose that $f_t$ is increasing on $[1,n]$ and take $x\in (n,n+1)$.
From (14) we have \\
$$\eqalign{
f_t'(x) & =t\frac {1}{x^2}-\frac {1}{x^2}\frac {1}{f_t(x-1)}+\left (1-\frac {1}{x}\right )
\frac {f_t'(x-1)}{f_t^2(x-1)} \cr
 & = \frac {1}{x^2f_t(x-1)}\left [tf_t(x-1)-1 \right ] +\left (1-\frac {1}{x}\right )
\frac {f_t'(x-1)}{f_t^2(x-1)} \cr
}$$
The second term is positive by the induction hypothesis, and the first term
is positive because \\ $tf_t(x-1)-1 \geq f_t(x-1)-1 \geq 0$, for all $x\geq 1$.
So, the function $f_t(x)$ is increasing on the interval $(n,n+1)$, and then, by
continuity, also on $[1,n+1]$. This completes the step of induction. \qed

As a consequence, we get the log-convexity for the sequence of {\bf central Delannoy
numbers}. Recall that the $n$-th central Delannoy number counts the number of lattice
paths in $(x,y)$ coordinate plane from $(0,0)$ to $(n,n)$ with steps $(1,0)$, 
$(0,1)$ and $(1,1)$. (Such paths are also known as {\bf king's paths}.)

{\bf Theorem 6}\\ (a) The sequence $D(n)$ of Delannoy numbers is log-convex.\\
(b) There exists $x=\lim_{n \rightarrow \infty} \frac {D(n)}{D(n-1)}$, and
$x=3+2 \sqrt 2$.

{\bf Proof}\\ (a) First note that the $n$-th central Delannoy number is the 
value of the $n$-th Legendre polynomial at $t=3$, $D(n)=P_n(3)$. This 
follows easily from the explicit expression for the 
generating function of the sequence $D(n)$, $D(x)=\frac {1}{\sqrt {1-6x+x^2}}$. 
Now apply Theorem 5.

(b) By (a) we know that $x_n(3)$ is increasing (and clearly bounded), and then
by passing to limit in (13) for $t=3$, the claim follows. \qed

%\bibliographystyle{plain}
%\bibliography{ref}

\end{document}